\documentclass[11pt,draft]{amsart}
\usepackage{amsmath,amssymb,amsfonts}

\newtheorem{conjecture}{Conjecture}[section]

\newtheorem{problem}{Problem}[section]

\theoremstyle{definition}
\newtheorem{definition}{Definition}[section]

\theoremstyle{remark}

\numberwithin{equation}{section}

\DeclareFontFamily{OT1}{wncyr}{\hyphenchar\font45 }
\DeclareFontShape{OT1}{wncyr}{m}{n}{%
  <5> <6> <7> <8> <9> gen * wncyr
   <10> <10.95> <12> <14.4> <17.28> <20.74>  <24.88>wncyr10}{}
\DeclareFontShape{OT1}{wncyr}{m}{it}{%
  <5> <6> <7> <8> <9> gen * wncyi
  <10> <10.95> <12> <14.4> <17.28> <20.74> <24.88> wncyi10}{}
\DeclareFontShape{OT1}{wncyr}{m}{sc}{%
  <5> <6> <7> <8> <9> <10> <10.95> <12> <14.4>
  <17.28> <20.74> <24.88>wncysc10}{}
\DeclareFontShape{OT1}{wncyr}{b}{n}{%
  <5> <6> <7> <8> <9> gen * wncyb
   <10> <10.95> <12> <14.4> <17.28> <20.74> <24.88>wncyb10}{}
\input cyracc.def
\def\rus{\usefont{OT1}{wncyr}{m}{n}\cyracc\fontsize{9}{11pt}\selectfont}

\DeclareMathSizes{9}{9}{7}{5} 

\begin{document}

\

\vskip -15mm

\centerline{\bf  Workshop ``Affine Algebraic
Geometry'',}

\centerline{\bf Oberwolfach, January 7--14, 2007}

\centerline{\bf Problems for Problem Session}

\vskip 10mm

\title[Quasihomogeneous affine threefolds]
{Birationally nonequivalent linear actions;\\Cayley
degrees of simple algebraic groups;\\
and singularities of two-dimensional quotients}

\author[Vladimir  L. Popov]{Vladimir  L. Popov${}^*$}
\address{Steklov Mathematical Institute,
Russian Academy of Sciences, Gubkina 8, 119991 Moscow,
Russia} \email{popovvl\char`\@orc.ru}

\thanks{
 ${}^*$\,Supported by Russian grants {\rus RFFI
05--01--00455}, {\rus N{SH}--9969.2006.1}, and program
{\it Contemporary Problems of Theoretical Mathematics}
of the Mathe\-matics Branch of the Russian Academy of
Sciences.}

\date{January 28, 2007}

\maketitle

\setcounter{tocdepth}{2}

 Below all algebraic varieties are taken over an
algebraically closed field $k$
of characteristic zero.
\section{Birationally nonequivalent linear actions}

Let $G$ be a reductive algebraic group. In 1992
 P. Katsylo published
 the following

 \begin{conjecture} {\rm(\cite{Ka})}\label{K1}
Let $V$ and $W$ be finite dimensional algebraic
$G$-modules with trivial stabilizers of points in
general position. Then the following properties are
equivalent:
\begin{enumerate}
\item[\rm(i)] $\dim V=\dim W$; \item[\rm(ii)] there
exists a $G$-equivariant birational map
$V\overset{\simeq}{\dashrightarrow} W$.
\end{enumerate}
\end{conjecture}

In \cite{Ka} Conjecture~\ref{K1} was proved for $G={\bf
SL}_2$, ${\bf PSL}_2$, and the symmetric groups ${\rm
S}_n$, $n\leqslant 4$. How\-ever E. Tevelev  observed
(unpublished) that Conjecture~\ref{K1} fails for
one-dimensional spaces and $G={\bf Z}/n$, $n\neq 2, 3,
4, 6$; the same observation was independently made in
\cite{RY}. In 2000 new counterexamples to
Conjecture~\ref{K1} have been found in \cite{RY}, where
a birational classification of finite dimensional
$G$-modules for diagonalizable $G$ has been obtained.
Being sceptical about Conjecture~\ref{K1},  in 1993 I
suggested to consider $W=V^*$, the dual module of $V$:

\begin{problem}\label{P1} Are there
a connected semisimple group $G$ and a finite
dimensional algebraic $G$-module $V$ with trivial
stabilizers of points in general position such that $V$
and $V^*$ are not birationally $G$-isomorphic?
\end{problem}

This problem was communicated to some people, see,
e.g.,\;\cite{RV}.
 So far it is still open.

It is well known that, if $G$ is connected, then $V^*$
is $V$ ``twisted'' by an automorphism of $G$. This
naturally leads to the following generalization. Let
$H$ be an algebraic group acting on an algebraic
variety $X$,
$$
H\times X\rightarrow X, \quad (h, x)\mapsto h\cdot x.
$$
Let $\sigma\colon H\rightarrow H$, $h\mapsto
{}^{\sigma}\!h$, be an automorphism of $H$. Consider
the following new action of $H$ on $X$:
$$
H\times X\rightarrow X, \quad (h, x)\mapsto
{}^{\sigma}\!h\cdot x.
$$
Then the new $H$-variety appearing in this way is
denoted by ${}^{\sigma}\!X$ and called $X$ ``twisted''
by $\sigma$. Problem \ref{P1} can be now generalized as
follows:
\begin{problem}\label{P1g} Are there
a connected semisimple group $G$, a finite dimensional
algebraic $G$-module $V$ with trivial stabilizers of
points in general position, and an automorphism
$\sigma$ of $G$ such that $V$ and ${}^{\sigma}V$ are
not birationally $G$-isomorphic?
\end{problem}

 Note that
if in Problem \ref{P1g} one replaces $V$  by a
$G$-variety $X$, then the answer is positive (see
\cite{RV} where this is proved for $G={\bf PGL}_n$).

 It is clear that  one should consider only outer
automorphisms $\sigma$ in Problem \ref{P1g}. Also, it
is easily seen that if $H$ is special in the sense of
J.-P. Serre, \cite{S},  then $X$ and ${}^{\sigma}\!X$
are always bi\-ra\-ti\-onally $G$-isomorphic,
cf.\;\cite{RV}. In particular, answering Problems
\ref{P1} and \ref{P1g} for simple $G$, one should
consider only
$${\bf SL}_n/\mu_d \mbox{ where }
d\neq 1, \mbox{ and the groups of types } {\sf D}_{l}
\mbox{ ($l$ is odd in Problem \ref{P1}), } {\sf E}_6.$$
In particular, we have the following special

\begin{problem} \label{Pg2} Let $G={\bf SL}_d/\mu_d$. Let
$V$ be the $d$-th symmetric power of $k^d$ endowed with
the natural action of $G$. Let $\sigma$ be the
automorphism of $G$ induced by the automorphism
$g\mapsto (g^{\sf T})^{-1}$ of $\;{\bf SL}_d$. Are $V$
and ${}^{\sigma}V$ birationally $G$-isomorphic?
\end{problem}

Note that  for $d>3$ in Problem \ref{Pg2}, stabilizers
of points in general position in $V$ are trivial,
\cite{P}.

\section{Cayley degrees of simple algebraic groups}

 Let $G$ be a connected reductive algebraic group
 and let ${\rm Lie} \,G$ be its Lie algebra.
 Consider the action of $G$ on ${\rm Lie} \,G$ via the adjoint representation
 and on $G$ by conjugation.
 \begin{definition}(\cite{LPR1}) $G$ is called {\it
 Cayley group} if $G$ and ${\rm Lie} \,G$ are
 birationally $G$-isomorphic.
 \end{definition}

  All simple Cayley groups have been
 classified in \cite[Theorem 1.31]{LPR1}: they are precisely the groups from the
 list
 $$
 {\bf SL}_n, n\leqslant 3;\hskip 3mm
 {\bf SO}_n, n\neq 2, 4;\hskip 3mm
 {\bf Sp}_{2n};\hskip 3mm
 {\bf PGL}_n.
 $$

 For every $G$, by  \cite[Prop. 10.5]{LPR1}
 there always exists a dominant
$G$-equivariant rational map $G\dashrightarrow {\rm
Lie} \,G$. So the following number is well defined:

\begin{definition}(\cite{LPR1}) The Cayley degree
${\rm Cay}(G)$ of $G$ is the minimum of degrees of
dominant rational $G$-equivariant maps
$G\dashrightarrow {\rm Lie} \,G$.
\end{definition}

So $G$ is Cayley if and only if ${\rm Cay}(G)=1$. In
general, ${\rm Cay}(G)$ ``measures'' how far $G$ is
from being Cayley.

\begin{problem}\label{lpr2} {\rm(\cite{LPR1})} Find the Cayley degrees of
connected simple algebraic groups.
\end{problem}

 In \cite{LPR1}, \cite{LPR2}) it is proved that
\begin{gather*}
{\rm Cay}({\bf SL}_n)\leqslant n-2,\hskip 1.5mm
\mbox{for $n\geqslant 3$};\quad
{\rm Cay}({\bf SL}_n/\mu_d)\leqslant n/d;\\
{\rm Cay} ({\bf Spin}_n)=
\begin{cases} 2&\text{for $n\geqslant 6$},\\
1&\text{for $n\leqslant 5$};
\end{cases}\\
{\rm Cay}({\bf G}_2)=2;\quad {\rm Cay}({\bf G}_2\times
{\bf G}_m^2)=1.
\end{gather*}

\noindent In particular, this implies that ${\bf
Cay}({\bf SL}_4)=2$ and $2\leqslant {\bf Cay}({\bf
SL}_5)\leqslant 3$.

\begin{problem} Find ${\rm Cay}({\bf SL}_5)$.
\end{problem}

At the moment no examples of groups whose Cayley degree
is bigger than $2$ are known.

\begin{problem} Is there $G$ such that ${\rm
Cay}(G)>2$? Is there a simple such $G$?
\end{problem}

\noindent More generally,

\begin{problem} Given a $d\in {\bf N}$, is there
$G$ such that ${\rm Cay}(G)>d$? Is there a simple
such~$G$?
\end{problem}

\section{Singularities of two-dimensional quotients}

Using a result of \cite{KR}, it was recently proved
 in \cite{G2} that
if a complex reductive algebraic group $G$
 acts algebraically on ${\bf C}^n$ and the categorical
quotient ${\bf C}^n/\!\!/G$ is two-dimensional, then
${\bf C}^n/\!\!/G$ is isomorphic to ${\bf C}^2/\Gamma$,
where $\Gamma$ is a finite group acting algebraically
on ${\bf C}^2$. This theorem can be considered as a
generalization of C.~T.~C. Wall's conjecture for the
linear action of $G$ on ${\bf C}^n$ proved
in~\cite{G1}.

This result, discussed in Koras' talk at this Workshop,
prompted the following question.

\begin{problem}{\rm (M. Miyanishi)} What are the groups
$\Gamma$ occuring in the above situation?
\end{problem}

I conjecture that  the following holds.

\begin{conjecture} If $G$ is connected,
then $\Gamma$ is cyclic.
\end{conjecture}

Note that it was conjectured in \cite{P1} and proved in
\cite{Ke} that if in the above situation the group $G$
is connected semisimple and the action of $G$ on ${\bf
C}^n$ is linear, then $\Gamma$ is trivial.


\begin{thebibliography}{LPRX}


\bibitem[G${}_1$]{G1}
R.~V.~Gurjar, {\it On a conjecture of C.~T.~C.~Wall},
J. Math. Kyoto Univ. {\bf 31} (1991), 1121--1124.


\bibitem[G${}_2$]{G2} R. Gurjar, {\em Two-dimensional
quotients of \;${\bf C}^n$ are isomorphic to ${\bf
C}^2/\Gamma$}, Transform. Groups {\bf 12} (2007), no.
1, 117--125.



\bibitem[Ka]{Ka} P. Katsylo, {\em On the birational
classification of linear representations}, preprint
MPI-92-1 (1992), Max-Planck-Inst. f\"ur Mathematik.

\bibitem[Ke]{Ke} G. Kempf, {\em Some quotient surfaces
are smooth} {\bf 27} (1980), 285--299.

\bibitem[KR]{KR} M.~Koras, P.~Russell, {\em
Contractible affine surfaces with a quotient
singularity}, Transform. Groups {\bf 12} (2007). no. 2,
to appear.

\bibitem[LPR${}_1$]{LPR1} N. Lemire, V. L. Popov,
Z.~Reichstein, {\em Cayley groups}, J. American
Mathema\-ti\-cal Society {\bf 19} (2006), 921--967.

\bibitem[LPR${}_2$]{LPR2} N. Lemire, V. L. Popov,
Z.~Reichstein, {\em On the Cayley degree of an
algebraic group}, in: {\em Actas del {\rm XVI} Coloquio
Latinoamericano de Algebra}, Revista Matematica
Iberoamericana, to appear, {\tt arXiv:} {\tt
math.AG/0608473}.

\bibitem[P]{P} A. M. Popov, {\em Irreducible simple linear
Lie groups with finite standard subgroups of general
position}, Funct. Anal. Appl. {\bf 9} (1975), 346--347.

\bibitem[P${}_1$]{P1} V. L. Popov, {\em Representations with
a free module of covariants}, Funct. Anal. Appl. {\bf
10} (1976), 242--244.




\bibitem[RV]{RV} Z. Reichstein, A. Vistoli, {\em Birational isomorphisms
between twisted group actions}, {\tt math.AG/}
{\tt
0504306.}

\bibitem[RY]{RY} Z. Reichstein, B. Youssin, {\em A birational invariant
for algebraic group actions}, Pacific J. Math. {\bf
204} (2002), no. 1, 223--246.

\bibitem[S]{S} J.-P. Serre, {\em Espaces fibr\'es
alg\'ebriques}, in: {\S\'eminaire Chevalley ``{\em
Anneaux de Chow}"}, expos\'e n${}^{\rm o}$\,1, ENS,
Paris, 1958, pp. 1-01--1-37.



\end{thebibliography}
\end{document}